\documentclass[12pt]{article} 
\def\R{\hbox{{\rm I}\kern-0.2em{\rm R}\kern0.2em}}%mathematical R for reals
\def\C{\hbox{{\rm I}\kern-0.2em{\rm C}\kern0.2em}}%mathematical C for complex
\def\D{\hbox{{\rm I}\kern-0.2em{\rm D}\kern0.2em}}

\def\be{\begin{equation}}
\def\ee{\end{equation}}

\def\({\left(}
\def\){\right)}
\def\[{\left[}
\def\]{\right]}
\def\bc{\begin{center}}
\def\ec{\end{center}}

\usepackage{epstopdf}
\usepackage{graphicx}
\usepackage{color}
\usepackage{float}
\def\C{\hbox{{\rm C}\kern-0.6em{\rm I}\kern0.6em}}

\begin{document}

\begin{center}{\large \bf Linearization from complex Lie point transformations}\end{center}

\bc{Sajid Ali$^{a}$, M. Safdar$^{b}$, Asghar Qadir$^{c}$\\

$^{a}$School of Electrical Engineering and Computer Science,
\\
$^{b}$School of Mechanical and Manufacturing Engineering, \\
$^{c}$School of Natural Sciences, \\
National University of Sciences and Technology, Campus H-12,
Islamabad 44000, Pakistan \\
sajid\_ali@mail.com, safdar.camp@gmail.com, aqadirmath@yahoo.com}\ec

\textbf{Abstract.} Complex Lie point transformations are used to
linearize a class of systems of second order ordinary differential
equations (ODEs) which have Lie algebras of maximum dimension $d$, with $d\leq 4$.
We identify such a class by employing complex
structure on the manifold that defines the geometry of differential
equations. Furthermore we provide a geometrical construction of the
procedure adopted that provides an analogue in $\R^{3}$ of the
linearizability criteria in $\R^2$.

\section{Introduction}
One method of solving a nonlinear ODE or a system of such equations is
to reduce it to linear form, that is called linearization, by
invertible transformations of the independent and dependent
variables (point transformations). There are many ways to convert to linear form like Lie point transformations, contact transformations, non-local transformations etc., but for convenience we will use the term ``linearization" to signify the Lie point transformations only and qualify the term if we mean any other transformation. Lie presented the most general
form of a scalar second order linearizable ODE by considering
arbitrary point transformations (see for example \cite{olver,ib} and the references therein). Over the
past few years there has been a rapidly growing interest in studying
linearization of higher-order ODEs and systems of these equations.
The simplest system ($f_{1}^{\prime\prime}=0, f_{2}^{\prime\prime}
=0$) is a unique system that admits a maximal algebra of dimension
$15=|sl(4,\R)|,$ which was proved in \cite{fels}, using Cartan's
equivalence method, after refining the result in \cite{gonzalez1}.
 A geometric approach
is used to obtain linearization criteria for systems of two second
order cubically semi-linear ODEs obtainable by projecting down a
system of three geodesic equations \cite{s5,aminova}. Indeed, utilizing
arbitrary point transformations, general forms of linearizable systems
and corresponding linearization criteria were studied in \cite{jm,
bagderina, meleshko}. The use of generalized Sundman transformations in the 
linearization problem is studied in \cite{dur}. The construction of
linearizing transformations from the first integrals of two
dimensional systems is carried out in \cite{chan, chan1}.

It is well known \cite{sb}, that to be linearizable by Lie point transformation a system of two
second order ODEs can have an algebra of maximum dimensions,
$5,6,7,8$ or $15$. The systems which do not have maximal algebras
are clearly non-linearizable in the sense that there exist no real
Lie point transformation which can be used to linearize them. This
paper addresses the problem of linearizability of such a class of
non-linearizable systems by complex Lie point transformations.

Recently we studied a special class of two dimensional linearizable
systems that corresponds to complex linearizable scalar ODEs
\cite{saf1, saj}. (The characterization for the correspondence will
be explained in the next section.) The linearizing transformations
to map nonlinear systems to linear forms are provided by the complex
fibre-preserving transformations
\begin{eqnarray}
\mathcal{L}_{1}:(x,u(x))\rightarrow (\chi(x),U(x,u)), \label{ftran}
\end{eqnarray}
where $u(x),$ is an analytic complex function of a real variable
$x$. The method was successfully applied to generate classes of
systems with maximum dimensions of their algebras $6,7$ and $15$, where the
linearizing transformations were obtained by the process of
``realification'' of the above transformations (\ref{ftran}). In this paper we
investigate a \emph{non-linearizable class} of systems of two ODEs
that can be obtained from a linearizable scalar complex equation,
but the complex linearizing transformations in this case are
different from those given above, i.e.,
\begin{eqnarray}
\mathcal{L}_{2}:(x,u(x))\rightarrow (\chi(x,u),U(x,u)), \label{ctran}
\end{eqnarray}
which is a complex Lie point transformation and contains (\ref{ftran}) as its special case. Notice that $\mathcal{L}_{2}$ cannot be used to obtain the
linearizing transformations for the corresponding systems. The
reason is that the transformed independent variable is complex which
gives two independent variables when splits into the real and
imaginary parts while the concerned systems are of ODEs rather than
PDEs. The standard techniques of linearization \cite{bagderina, meleshko,sb}, can not be
applied to such systems as they have symmetry algebra of maximum dimension
$d\leq 4$. However we prove that linearization can still be achieved
for such systems in the complex plane and the solutions can also be
obtained explicitly. Most surprisingly, we find that there can be systems that have no Lie point symmetries but can be solved by this procedure of converting to a complex scalar ODE and linearizing the scalar equation to write down the solution in the complex domain and thence obtain the solution of the system. An example of such a system is also given. We call this procedure for linearizing systems of
two second order ODEs from complex point transformations with less than
five symmetries \emph{complex linearization}.

The procedure of constructing a pair of real functions of two
variables from a single complex function of a complex variable,
leading to a system of PDEs, entails the use of the Cauchy-Riemann
(CR) equations in a transparent way. However, the role of these
equations for a system of ODEs is far from clear as the CR-equations
require two independent variables. Here we give a simple explanation
of this role using geometry, for the free particle
ODE which presents an elegant framework of
viewing the straight line (which is maximally symmetric and
invariant under the full group $SL(3,\R)$) in a higher dimensional
space which is only possible if we put on \emph{complex
glasses}. The CR-equations are shown to play
an essential part in establishing the correspondence between
solutions of the base complex linearizable or integrable ODEs and
emerging systems.

The outline of the paper is as follows. The criteria for the
correspondence of systems with complex linearizable equations are
given in the second section. The third section is on the
CR-equations associated with systems of two second order ODEs and
exploring their role in establishing correspondence of solutions of
systems and base scalar ODEs. The subsequent section contains
application of complex linearization procedure on systems which
correspond to linearizable complex base ODEs. The last section is
devoted to the conclusion and discussion.

\section{Complexification of Systems of ODEs and Classification}

We first explain the basic formalism of complex linearizability by
taking a general system of two second order ODEs
\begin{eqnarray}
f_{1}^{\prime\prime}=\omega_{1}(x,f_{1},f_{2},f_{1}^{\prime},f_{2}^{\prime}),\nonumber\\
f_{2}^{\prime\prime}=\omega_{2}(x,f_{1},f_{2},f_{1}^{\prime},f_{2}^{\prime}),
\label{odesys}
\end{eqnarray}
which may be regarded as a surface $S$ in a $2(3)+1=7-$dimensional space
whose components comprising of the independent and dependent
variables along with their derivatives. Hence a solution of system
(\ref{odesys}) is an integral curve on surface $S$. We now introduce
a complex structure $J_{_{\C}}:\R^{6}\rightarrow \C^{3},$ on the
$6-$dimensional subspace of $S$ by assuming
$f_{1}(x)+if_{2}(x)=u(x),$ where all first order and second order
derivatives of the real functions $f_{1},~f_{2}$ of a real variable
$x$ are determined with $u^{\prime}$ and $u^{\prime\prime}$,
respectively. Therefore our solution curve is now embedded in a
complex $3_{_{\C}}-$dimensional space, $\C^3$. If we regard $x$ as
the dimension of time then this can be viewed as the propagation of
our solution curve in time in a complex space $\C^3$. This
identification gives deeper insights into the symmetry analysis
which we shall see as we proceed to the subsequent sections. This
yields a class of those two dimensional systems (\ref{odesys}) which
can be projected to a scalar second order complex equation
\begin{eqnarray}
 u^{\prime\prime}=\omega(x,u,u^{\prime}), \label{code}
\end{eqnarray}
where $\omega(x,u)=\omega_{1}(x,u)+i\omega_{2}(x,u),$ in a
$3_{_{\C}}+1=4-$dimensional partially complex space $S_{_{\C}}$,
which is comprised of a $3_{_{\C}}-$complex dimensional subspace and
a one-dimensional subspace that correspond to independent variable
$x$ which is not complex. The basic criterion to identify such
systems is that both $\omega_{1}$ and $\omega_{2}$ in
(\ref{odesys}), satisfy the CR-equations
\begin{eqnarray}
\omega_{1,f_{1}}=\omega_{2,f_{2}}~,~~\omega_{1,f_{2}}=-\omega_{2,f_{1}}, \nonumber \\
\omega_{1,f_{1}^{\prime}}=\omega_{2,f_{2}^{\prime}}~,~~\omega_{1,f_{2}^{\prime}}=-\omega_{2,f_{1}^{\prime}},
\label{cr0}
\end{eqnarray}
namely, both functions are analytic of their arguments. This line of approach has been
followed in \cite{saf1,saj} to characterize systems of ODEs that
emerge from linearizable complex equations. Several non-trivial and
interesting results are obtained for two dimensional systems
\emph{despite} ``trivial'' identification of systems (\ref{odesys})
from scalar equations (\ref{code}).

A two dimensional system of ODEs
\begin{eqnarray}
f_{1}^{^{\prime \prime }}=A_{1}f_{1}^{^{\prime
}3}-3A_{2}f_{1}^{^{\prime }2}f_{2}^{^{\prime
}}-3A_{1}f_{1}^{^{\prime }}f_{2}^{^{\prime }2}+A_{2}f_{2}^{^{\prime
}3}+B_{1}f_{1}^{^{\prime }2}-2B_{2}f_{1}^{^{\prime }}f_{2}^{^{\prime
}}-B_{1}f_{2}^{^{\prime }2}+\nonumber\\
C_{1}f_{1}^{^{\prime
}}-C_{2}f_{2}^{^{\prime }}+D_{1},\nonumber\\
f_{2}^{^{\prime \prime }}=A_{2}f_{1}^{^{\prime
}3}+3A_{1}f_{1}^{^{\prime }2}f_{2}^{^{\prime
}}-3A_{2}f_{1}^{^{\prime }}f_{2}^{^{\prime }2}-A_{1}f_{2}^{^{\prime
}3}+B_{2}f_{1}^{^{\prime }2}+2B_{1}f_{1}^{^{\prime }}f_{2}^{^{\prime
}}-B_{2}f_{2}^{^{\prime }2}+\nonumber\\
C_{2}f_{1}^{^{\prime
}}+C_{1}f_{2}^{^{\prime }}+D_{2}, \label{cubicsys}
\end{eqnarray}
where $A_j,~B_j,~C_j,~D_j~ (j=1,2)$ are functions of $x$ and $f_j$,
is a candidate of complex linearization if and only if its
coefficients satisfy a set of four constraint equations given in
\cite{saj}. It was shown that the system (\ref{cubicsys})
corresponds to a cubically semi-linear complex scalar equation of
the form
\begin{eqnarray}
u^{\prime\prime}+E_{3}(x,u)u^{\prime 3}+E_{2}(x,u)u^{\prime
2}+E_{1}(x,u)u^{\prime}+E_{0}(x,u)=0. \label{cubiccode}
\end{eqnarray}
where $E_k,~(k=0,1,2,3)$ are complex functions of $x$ and $u$, which
according to Lie is linearizable using complex point transformations provided the coefficients
satisfy linearization criteria given in \cite{ib}. Therefore system
(\ref{cubicsys}) is complex linearizable if and only if the associated base
equation (\ref{cubiccode}) can be mapped to $U^{\prime\prime} =0$,
using a complex Lie point transformation (\ref{ctran}). It may be pointed out that there is only one candidate of complex linear equation which is $U^{\prime\prime} =0$, unlike systems of two linear ODEs which have five candidates corresponding to their algebras $5,6,7,8$ or $15$ of maximum dimensions. Therefore to be linearizable by real Lie point transformations system (\ref{cubicsys}) can have one of these algebras.

The systems of the form (\ref{cubicsys}) that correspond to a complex linear equation using fibre-preserving transformations (\ref{ftran}) were studied in \cite{saf1}, where it was found that they are at most quadratic in the first derivatives and can have algebras of maximum dimensions 6, 7 or 15. They can also be linearized using real Lie point transformations by standard techniques because the dimensions of their algebras coincide with the dimensions of linearizable systems. We denote this class of quadratically semi-linear systems by $\Upsilon_{1}$. It turns out that if we take systems (\ref{cubicsys}) with cubic dependence on the first derivatives they are no longer linearizable because the dimensions
of their Lie algebras do not coincide with the dimensions of linearizable systems. Therefore we obtain a class $\Upsilon_{2}$, of systems of the form (\ref{cubicsys}), with cubic nonlinearity and such that the maximum dimensions of their Lie algebras can be 4, 3, 2, 1 or 0. Such systems are not \emph{not} linearizable but they correspond to a linearizable (complex) scalar second order ODE. Now all systems of the form (\ref{cubicsys}) for which coefficients satisfy constraint equations \cite{saj}, correspond to complex linear equation (\ref{cubiccode}) under the complex point transformation (\ref{ctran}), therefore all systems in $\Upsilon_{2}$, are also complex linearizable. Since the dimensions of Lie algebras of systems in classes $\Upsilon_1$ and
$\Upsilon_2$ are not equal, therefore there do not exist any complex point
transformations that map a system in $\Upsilon_1$ to a system in
$\Upsilon_2$ and vice versa, therefore
\begin{eqnarray}
\Upsilon_1 \cap \Upsilon_2 = \emptyset.
\end{eqnarray}
In short, the class $\Upsilon_1$ contains those complex linearizable
systems which are at most quadratically semi-linear and they can be
mapped to the optimal canonical form \cite{saf1}, whereas all
complex linearizable systems with cubic nonlinearity are contained
in $\Upsilon_2$. Below we summarize two classes arising from a complex linear equation in Table 1.
\begin{table}[H]
\caption{Classification of two-dimensional systems arising from linear complex equation}
\begin{center}
\begin{tabular}{|c|c|c|c|c|}
 \hline
  Class & Lie Algebra Dimensions & Linearizable & Complex Linearizable & Type  \\ \hline \hline
  $\Upsilon_{1}$ & 6,7,15 & Yes & Yes & Quadratic \\ \hline
  $\Upsilon_{2}$ & 0,1,2,3,4 & No & Yes & Cubic \\
  \hline
\end{tabular}
\end{center}
\end{table}
Note that the cases of Lie algebras of maximum dimensions 5 and 8 are not contained in $\Upsilon_1\bigcup\Upsilon_2$ but do give linearizable systems. For completeness we also discuss a class of those systems that can be solved from complex methods in which systems are solved due to their correspondence with complex solvable equation. It is given in the Appendix.

In order to construct the representative system for the class
$\Upsilon_2$ with fewer symmetries \emph{yet} is complex
linearizable we focus on a class of cubically semi-linear systems
 \begin{eqnarray}
f_{1}^{\prime\prime}=\beta f_{1}^{\prime 3}-3\gamma f_{1}^{\prime
2}f_{2}^{\prime}-3\beta f_{1}^{\prime}f_{2}^{\prime
2}+\gamma f_{2}^{\prime 3},\nonumber\\
f_{2}^{\prime\prime}=\gamma f_{1}^{\prime 3}+3\beta f_{1}^{\prime
2}f_{2}^{\prime}-3\gamma f_{1}^{\prime}f_{2}^{\prime 2}-\beta
f_{2}^{\prime 3}, \label{cubic}
\end{eqnarray}
from (\ref{cubicsys}), where $\beta=\beta(x,f_{1},f_{2})$ and
$\gamma=\gamma(x,f_{1},f_{2})$. The complex linearizability criteria \cite{saj} is satisfied for
these systems if and only if the coefficients satisfy
\begin{eqnarray}
\beta_{xx}=0, \quad \gamma_{xx}=0,
\end{eqnarray}
i.e., $\beta=b_1x+b_2,~ \gamma=c_1x+c_2,$ therefore the system (\ref{cubic}) is complex
linearizable. It is easy to verify that the systems (\ref{cubic})
have symmetry algebras of maximum dimensions less than $4$, provided all constants $b_1,b_2,c_1$ and $c_2$ do not vanish
simultaneously. The system (\ref{cubic})
satisfy (\ref{cr0}) and thus can be projected to (\ref{cubiccode}),
\begin{eqnarray}
u^{\prime\prime}+E_{3}(x,u)u^{\prime3}=0,
\end{eqnarray}
with $E_2=0=E_1=E_0,$ which can be linearized to $U^{\prime\prime}=0$, using a complex point transformation
(\ref{ctran}). Before proceeding to the
applications and characterizing such systems we first describe the geometry of $U^{\prime\prime}=0$, under
general complex Lie point transformations (\ref{ctran}).

\section{Geometry of Complex Linearization}

We know that all linearizable scalar differential equations are
equivalent to the free particle equations (see, e.g., \cite{olver,ib})
whose solution is a straight line. The crucial step after ensuring
complex linearizability is to obtain the transformations which help
in the integration of the systems in $\Upsilon_1$ and $\Upsilon_2$.
In \cite{saf1}, the complex fibre-preserving transformations (\ref{ftran}),
were used to map a system in class $\Upsilon_1$ into the free
particle complex equation, $U^{\prime\prime}=0,$ where prime denotes
differentiation with respect to $\chi$. The real and imaginary parts
of such a free particle equation yields a system
\begin{eqnarray}
F_{1}^{\prime\prime}=0,~~F_{2}^{\prime\prime}=0,\label{linsys}
\end{eqnarray}
where $U=F_{1} + i F_{2}$, therefore the corresponding system can be mapped to system (\ref{linsys}), where the
complex transformations (\ref{ftran}) were used to obtain the real linearizing transformations. It is noteworthy that the transformations in this case are
$(real,complex)\rightarrow(real,complex)$. On the other hand a general
complex point transformation is of the form
\begin{eqnarray}
\mathcal{L}:(x,u)\rightarrow (\chi(x,u),U(x,u)),\nonumber
\end{eqnarray}
namely, $(real,complex)\rightarrow(complex,complex)$, in which the
first argument $\chi,$ can be a complex function thereby adds a
superficial dimension. It makes a huge difference
on linearizability of systems in $\Upsilon_{2}$. In particular, upon splitting it into the
real and imaginary parts
\begin{eqnarray}
\chi(x,u)=\chi_1(x,f_1,f_2)+i\chi_2(x,f_1,f_2), \label{chi}
\end{eqnarray}
and since the dependent function $U(\chi)$ is complex which yields
two real functions $F_1$ and $F_2,$ both of them are not only
functions of $\chi_1$ but also of $\chi_2$, i.e.,
\begin{eqnarray*}
U(\chi)=F_{1}(\chi_1,\chi_2) + i F_{2}(\chi_1,\chi_2), \label{chi}
\end{eqnarray*}
therefore the linearized
scalar equation, $U^{\prime\prime}=0,$ fails to produce the free
particle system (\ref{linsys}). Notwithstanding, the solution of a
system in class $\Upsilon_2$ is
extractable from the complex solution $U(x)$, upon its split. We now give a geometrical understanding
of complex linearization associated with base equation $U^{\prime\prime} = 0$. Since the prime denotes differentiation with
respect to $\chi$ which upon using the chain rule yields
\begin{eqnarray}
\frac{\partial}{\partial\chi}=
\frac{1}{2}\left(\frac{\partial}{\partial\chi_{_1}}-i\frac{\partial}{\partial\chi_{_2}}\right),
\end{eqnarray}
and so $U^{\prime\prime}=0$, is a system of two partial differential equations
\begin{eqnarray}
F_{_{1_{\chi_1\chi_1}}}-F_{_{1_{\chi_2\chi_2}}}+2F_{_{2_{\chi_1\chi_2}}}=0,\nonumber\\
F_{_{2_{\chi_1\chi_1}}}-F_{_{2_{\chi_2\chi_2}}}-2F_{_{1_{\chi_1\chi_2}}}=0.
\label{linpde}
\end{eqnarray}
Now by definition a complex Lie point transformation is analytic
thus $\mathcal{L}$ is analytic. Since the derivative $u^\prime$
transforms into a complex derivative $U^\prime$ which exists if and
only if $U(\chi)$ is complex analytic and is preserved under
$\mathcal{L}$ therefore
\begin{eqnarray}
F_{_{1_{\chi_1}}}=F_{_{2_{\chi_2}}},\quad
F_{_{1_{\chi_2}}}=-F_{_{2_{\chi_1}}},\label{cr}
\end{eqnarray}
which are the CR equations. The solution of system (\ref{linpde})
along with the condition (\ref{cr}) upon using invertible
transformations (\ref{ctran}) reveal solutions of the original
system. Hence we have established the following result.
\newline
\newline
\textbf{Theorem 1.} \emph{All complex linearizable two dimensional
systems of ODEs in class $\Upsilon_2$ can be transformed
into systems of PDEs} (\ref{linpde}) \& (\ref{cr}) \emph{under the
transformation} (\ref{ctran}).
\newline
\newline
We now develop some geometrical aspects
of complex linearization. In order to do that we employ the original
idea of Riemann, that a complex mapping may be regarded as the
dependence of one plane on another plane, unlike the dependence of a
real function on a line. To do this we first obtain the solution of
system (\ref{linpde})$-$(\ref{cr}) which is
\begin{eqnarray}
F_{1}(\chi_1 , \chi_2 )=c_1\chi_1 +c_2\chi_2 +c_3, \nonumber\\
F_{2}(\chi_1 , \chi_2 ) = c_1\chi_2 -c_2\chi_1 +c_4, \label{plane}
\end{eqnarray}
where $c_m,~ (m=1,2,3,4)$ are real arbitrary constants. These are
two coordinate planes determined by $\chi_1$ and $\chi_2$ with
normals
\begin{eqnarray}
\mathbf{n}_1 = [c_1,~ c_2] , \nonumber\\
\mathbf{n}_2 = [c_2, -c_1] ,
\end{eqnarray}
thus they intersect at right angles
\begin{eqnarray}
\mathbf{n}_1 \cdot \mathbf{n}_2 = 0 ,
\end{eqnarray}
resulting in a straight line at intersection. Thus the geometric
linearizing criterion for scalar second order differential
equations, namely a straight line, is extended to the intersection
of two planes at right angle in the complex linearization of two
dimensional systems. Note that both $\chi_1$ and $\chi_2$ in
(\ref{chi}) are functions of $(x,f_1,f_2).$ Therefore the role of
$\chi_2$ can be regarded as slicing the three-dimensional space
$\R^3=\{(x,f_1,f_2)\}$ into two coordinate planes. Interestingly, the
solution $(f_1,f_2)$ of the system under consideration is found by
solving (\ref{plane}) with the use of $F_1$ and $F_2$ from $U(\chi)$
in (\ref{ctran}). Hence we arrive at the following geometrical
result.
\newline\newline \newline
\textbf{Theorem 2.} \emph{The necessary and sufficient condition for
a two dimensional system} (\ref{cubicsys}) \emph{to be
complex linearizable is that the two planes determined by}
(\ref{plane}) \emph{intersect at right angle resulting in a straight
line which corresponds to scalar linear equations.}

Figure 1 illustrates the geometry and presents an elegant
description of complex linearization.
\begin{figure}[H]
\centerline{\includegraphics[width=15cm]{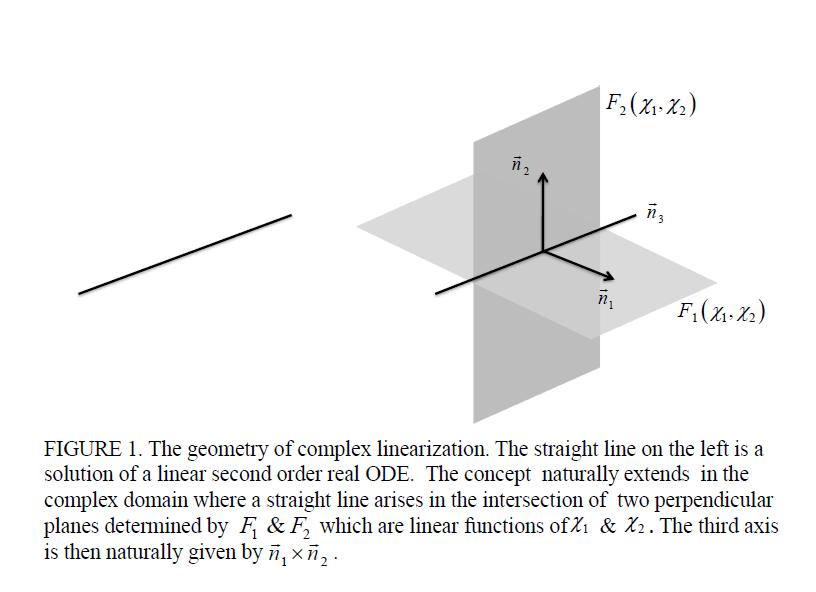}}
\end{figure}
 \textbf{Illustrative Example:}\quad We
consider an example of a physical system known as coupled modified Emden system \cite{ib},
\begin{eqnarray}
f_1^{\prime\prime}=-3f_1f_1^{\prime}+3f_2f_2^{\prime}-f_1^3+3f_1f_2^2, \nonumber \\
f_2^{\prime\prime}=-3f_2f_1^{\prime}-3f_1f_2^{\prime}+f_2^3-3f_1^2f_2.
\label{sys7}
\end{eqnarray}
to explain briefly how the procedure works. This system has three symmetries {$X_1,~X_2,~X_3$}, where
\begin{eqnarray}
X_{1}=\frac{\partial}{\partial{x}},\quad
X_{2}=x\frac{\partial}{\partial{x}}-
f_1\frac{\partial}{\partial{f_{1}}}-f_2\frac{\partial}{\partial{f_{2}}}, \nonumber\\
X_{3}=x^2\frac{\partial}{\partial{x}}-2x
f_1\frac{\partial}{\partial{f_{1}}}-2xf_2\frac{\partial}{\partial{f_{2}}},
\end{eqnarray}
with Lie algebra $[X_1,X_2]=X_1$, $[X_1,X_3]=2X_2$ and
$[X_2,X_3]=X_3$. The system (\ref{sys7}) is solvable \emph{only} by
complex linearization. The complex point transformation
\begin{equation}
\chi =x-\frac{1}{u},\quad
U=\frac{x^{2}}{2}-\frac{x}{u},\label{tran0}
\end{equation}
is of the form (\ref{ctran}), which does the trick to map the associated complex equation (\ref{code}) of (\ref{sys7}), into the complex free particle
equation, whose solution after inverting the above transformations
directly yields the solution of the system (\ref{sys7}).
In this case (\ref{linpde}) and (\ref{cr}) yield
\begin{eqnarray}
F_1 = a_1\chi_1 -a_2\chi_2 +b_1,\nonumber\\
F_2 = a_2\chi_1 +a_1\chi_2 +b_2, \label{tran1}
\end{eqnarray}
where
\begin{eqnarray}
&\chi_1 =x-\frac{f_1}{f_1^2+f_2^2}~, \quad
\chi_2=\frac{f_2}{f_1^2+f_2^2}~,\nonumber
\\ &F_1 = \frac{x^2}{2}-\frac{xf_1}{f_1^2+f_2^2}~,\quad F_2 =
\frac{xf_2}{f_1^2+f_2^2}~, \label{tran3}
\end{eqnarray}
obtained from (\ref{tran0}). Now by solving (\ref{tran1}), for $f_1$
and $f_2$, and by invoking equations (\ref{tran3}), we get the
solution
\begin{eqnarray}
f_1(x)=\frac{2x^3-6x^2a_1+4(a_2^2+a_1^2-b_1)x+
4a_1b_1+4a_2b_2}{x^4-4x^3a_1+4((a_2^2+
a_1^2-b_1)x^2+2(a_2b_2+a_1b_1)x+b_1^2+b_2^2)}, \nonumber\\
f_2(x)=\frac{(2x^2+4b_1)a_2+4b_2(x-a_1)}{x^4-4x^3a_1+4((a_2^2+
a_1^2-b_1)x^2+2(a_2b_2+a_1b_1)x+b_1^2+b_2^2)},
\end{eqnarray}
of system (\ref{sys7}).

\section{Applications}

We now illustrate the theory with the aid of examples on
class $\Upsilon_2$, i.e., the
complex linearizable systems with algebras of dimensions
$d\leq 4$. It may be pointed out that systems which can not
be dealt via standard Lie symmetry approach
are solvable via symmetry approach! We also highlight the procedure
with which we obtained these systems systematically.
\newline
\newline
\textbf{a. Solvable system of $4-$dimensional algebra:}\quad
Considering $\beta(x,f_{1},f_{2})=1$ and $\gamma(x,f_{1},f_{2})=0,$
in (\ref{cubic}) we obtain a coupled system
\begin{eqnarray}
f_{1}^{\prime\prime}-f_{1}^{\prime 3}+3f_{1}^{\prime}f_{2}^{\prime
2}=0,\nonumber\\
f_{2}^{\prime\prime}-3f_{1}^{\prime 2}f_{2}^{\prime}+f_{2}^{\prime
3}=0, \label{sys3}
\end{eqnarray}
which is complex linearizable and has only four symmetries
\begin{eqnarray}
X_{1}=\frac{\partial}{\partial{x}},\quad
X_{2}=\frac{\partial}{\partial{f_{1}}},\quad X_{3}=
\frac{\partial}{\partial{f_{2}}},\quad
X_{4}=2x\frac{\partial}{\partial{x}}+
f_{1}\frac{\partial}{\partial{f_{1}}}+f_{2}\frac{\partial}{\partial{f_{2}}},
\end{eqnarray}
with Lie algebra
\begin{eqnarray}
~[X_1,X_2]=0,\quad [X_1,X_3]=0,\quad [X_2,X_3]=0,\nonumber\\
~[X_1,X_4]=2X_1,\quad [X_2,X_4]=X_2,\quad [X_3,X_4]=X_3,
\end{eqnarray}
therefore, it is not in one of the linearizable classes of two
dimensional systems. Now in order to carry out integration of system
(\ref{sys3}), we linearize the corresponding equation
\begin{eqnarray}
u^{\prime\prime}-u^{\prime 3}=0,
\end{eqnarray}
which has an $8-$dimensional Lie algebra. This can be mapped to the
linear equation
\begin{eqnarray}
U^{\prime\prime}+1=0,
\end{eqnarray}
by inverting the role of the independent and dependent variables
$\chi=u,~U=x$. It has the solution $2U=-\chi^2+a\chi+b,$ where $a$
and $b$ are complex constants, which in terms of the original
variables becomes $u(x)=\pm\sqrt{a-2x}+b,$ and yields the solution
\begin{eqnarray}
f_1(x)=\pm\left (\frac{~a_1-2x+\sqrt{(a_1-2x)^2+a_2^2}}{2}\right )^{1/2}+b_1,      \nonumber\\
f_2(x)=\pm\left (\frac{-a_1+2x+\sqrt{(a_1-2x)^2+a_2^2}}{2}\right
)^{1/2}+b_2,
\end{eqnarray}
of (\ref{sys3}). System (\ref{sys3}) can be solved by real
symmetry analysis because it has four symmetry generators. However,
we now proceed to systems that are {\it not} solvable by real
symmetry methods as they have \emph{fewer symmetry generators than
four}.
\newline
\newline
\textbf{b. Solvable system of $3-$dimensional algebra:}\quad It is
easy to construct a system from (\ref{cubic}) which has only three
symmetries. For example we observe that in (\ref{cubic}) the
functions $\beta$ and $\gamma$ can be at most linear functions of
independent variable $x$. Hence we obtain a complex linearizable
system
\begin{eqnarray}
f_{1}^{\prime\prime}-xf_{1}^{\prime 3}+3xf_{1}^{\prime}f_{2}^{\prime
2}=0,\nonumber\\
f_{2}^{\prime\prime}-3xf_{1}^{\prime 2}f_{2}^{\prime}+xf_{2}^{\prime
3}=0,\label{sys4}
\end{eqnarray}
by involving $x$ linearly in the coefficients to remove the
$x-$translation. Thus we obtain the following $3-$dimensional
Abelian Lie algebra
\begin{eqnarray}
X_{1}=x\frac{\partial}{\partial{x}},\quad
X_{2}=\frac{\partial}{\partial{f_{1}}},\quad X_{3}=
\frac{\partial}{\partial{f_{2}}}.
\end{eqnarray}
We follow the same procedure as developed in the previous case and
solve the corresponding nonlinear equation
\begin{eqnarray}
u^{\prime\prime}-xu^{\prime 3}=0.
\end{eqnarray}
It is transformable to a linear form $U^{\prime\prime}=-U,$ which
after inverting the variables gives
\begin{eqnarray}
u(x)=\arctan\left( \frac{x}{\sqrt{a-x^2}} \right) +b.
\end{eqnarray}
System (\ref{sys4}) is {\it not} solvable by real symmetry methods
but by complex linearization.
\newline
\newline
\textbf{c. Solvable system of $2-$dimensional algebras:}\quad
Consider the system
\begin{eqnarray}
f_{1}^{\prime\prime}-f_{1}f_{1}^{\prime 3}+3f_{2}f_{1}^{\prime
2}f_{2}^{\prime}+3f_{1}f_{1}^{\prime}f_{2}^{\prime
2}-f_{2}f_{2}^{\prime 3}=0,\nonumber\\
f_{2}^{\prime\prime}-f_{2}f_{1}^{\prime 3}-3f_{1}f_{1}^{\prime
2}f_{2}^{\prime}+3f_{2}f_{1}^{\prime}f_{2}^{\prime
2}+f_{1}f_{2}^{\prime 3}=0, \label{sys5}
\end{eqnarray}
which has only two Lie symmetries
\begin{eqnarray}
X_{1}=\frac{\partial}{\partial{x}},\quad
X_{2}=3x\frac{\partial}{\partial{x}}+
f_1\frac{\partial}{\partial{f_{1}}}+f_2\frac{\partial}{\partial{f_{2}}}.
\end{eqnarray}
The system (\ref{sys5}) is solvable due to its correspondence with
the complex scalar second order ODE
\begin{eqnarray}
u^{\prime\prime}-uu^{\prime 3}=0,
\end{eqnarray}
which {\it is} linearizable, despite having only a two dimensional
algebra. Notice that even a scalar second order ODE requires at
least two symmetries to be solvable, while here complex
linearization helps to solve a system of two ODEs. In the subsequent
cases we provide the solution for a system with {\it only one}
symmetry generator, {\it which is insufficient to solve even a
scalar second order ODE}. Nevertheless, {\it we can go further} !
\newline
\newline
\textbf{d. Solvable system of $1-$dimensional algebra:}\quad
Consider $\beta(x,f_{1},f_{2})=xf_{1}$ and
$\gamma(x,f_{1},f_{2})=xf_{2},$ in (\ref{cubic}) we obtain
\begin{eqnarray}
f_{1}^{\prime\prime}-xf_{1}f_{1}^{\prime 3}+3xf_{2}f_{1}^{\prime
2}f_{2}^{\prime}+3xf_{1}f_{1}^{\prime}f_{2}^{\prime
2}-xf_{2}f_{2}^{\prime 3}=0,\nonumber\\
f_{2}^{\prime\prime}-xf_{2}f_{1}^{\prime 3}-3xf_{1}f_{1}^{\prime
2}f_{2}^{\prime}+3xf_{2}f_{1}^{\prime}f_{2}^{\prime
2}+xf_{1}f_{2}^{\prime 3}=0.\label{sys6}
\end{eqnarray}
This system is non-linearizable as it has \emph{only} a scaling
symmetry $X_{1}=x\partial_{x}$. The corresponding scalar second
order complex ODE is
\begin{eqnarray}
u^{\prime\prime}-xuu^{\prime 3}=0,
\end{eqnarray}
which has an $8-$dimensional algebra and linearizes to
\begin{eqnarray}
U^{\prime\prime}+\chi U=0,
\end{eqnarray}
which is the Airy equation whose solutions are Airy functions
extended to the complex plane. The solution of the complex
linearized equation for $U(x)$ is given by
\begin{eqnarray}
U(\chi)=c_{1}\mbox{Ai}(-\chi)+c_{2}\mbox{Bi}(-\chi),\label{sol}
\end{eqnarray}
where $\mbox{Ai}(-\chi)$ and $\mbox{Bi}(-\chi),$ are the two Airy
functions. Inverting ($\ref{sol}$), we obtain a solution of the
associated nonlinear equation which implicitly provides a solution
\begin{eqnarray}
\Re(c_{1}\mbox{Ai}(-f_{1}-if_{2})+c_{2}\mbox{Bi}(-f_{1}-if_{2}))=x,\nonumber\\
\Im(c_{1}\mbox{Ai}(-f_{1}-if_{2})+c_{2}\mbox{Bi}(-f_{1}-if_{2}))=0,
\end{eqnarray}
where $\Re$ and $\Im$ are the real and imaginary parts of the
arguments, for the system ($\ref{sys6}$).

\textbf{e. Solvable system of $0-$dimensional algebra:}\quad
We now take the quadratic freedom of $f_{1}$ and $f_{2}$, in $\beta$ and $\gamma$, to get rid of
the remaining symmetry in the above system. Consider $\beta=x-f_{1}^2 +f_{2}^2$ and
$\gamma=2f_{1}f_{2},$ in (\ref{cubic}), we obtain
\begin{eqnarray}
f_{1}^{\prime\prime}=(x^2-f_{1}^2+f_{2}^2) (f_{1}^{\prime 3}-3f_{1}^{\prime}f_{2}^{\prime 2} )+2 f_{1}f_{2}(3 f_{1}^{\prime
2}f_{2}^{\prime} -f_{2}^{\prime 3}),\nonumber\\
f_{2}^{\prime\prime}=(x^2-f_{1}^2+f_{2}^2) (3f_{1}^{\prime 2}
f_{2}^{\prime }-f_{2}^{\prime 3})- 2 f_{1}f_{2}(f_{1}^{\prime 3}-3
f_{1}^{\prime }f_{2}^{\prime 2} ), \label{sys8}
\end{eqnarray}
This system has no real point symmetry. However the corresponding complex
equation
\begin{eqnarray}
u^{\prime\prime}-xu^{2}u^{\prime 3}=0,
\end{eqnarray}
is again linearizable to $U^{\prime\prime}=0$, producing the
solution of system (\ref{sys8}).

\section*{Conclusion}
Complex symmetry analysis provides a class of systems of two ODEs obtainable 
from a scalar second order equation if the dependent variable is a
complex function of a real independent variable. The linearizability
of this base complex equation generates two different classes of
systems of two second order ODEs: (a) real linearizable systems that can be linearized from real Lie point transformations; and
(b) complex linearizable systems that can only be linearized using complex Lie point transformations. In \cite{saf1}, it was shown that systems that can be transformed into a complex linear equation has $6$, $7$ or $15-$ maximum dimensional algebra which can, therefore, also be linearized by real Lie point transformations. The second class is investigated here and it is found that these systems
are not linearizable by any real Lie point transformations but they are solvable due to their correspondence
with linearizable scalar complex equations. This class contains
complex linearizable systems with maximum dimensions of Lie algebras less than four.
The complex point transformations play a significant role in their
linearization.

It would be worthwhile to investigate the contact symmetries of given systems
that might provide deeper insights into complex point transformations. The group of
contact transformations is infinite-dimensional for linear systems and contains the
group of fibre-preserving transformations and point transformations as its special cases. The
last system with no symmetry clearly indicate an inherited beauty
in the complex domain which is not visible in the real domain. It is a benchmark problem as to why this system is solvable in the complex
domain and not in the real.

Important features of the given systems can be explored in terms of first integrals by applying the Noether symmetry analysis, provided the given systems come from a variational principle. In contrast, Cartan's equivalence approach (which does not necessarily depend on the existence of a Lagrangian) can be applied to investigate invariants and differential invariants of such systems. Further directions may include the problem of linearization of even dimensional systems of ODEs from complex point transformations.

\section*{Appendix}

For completeness we also examine the solvability of those systems that
are not complex linearizable yet can be solved via complex procedure
if they are mapped to solvable (integrable) scalar complex
equations with two dimensional solvable algebras. For this purpose
we state the general form of scalar equations with two symmetries.
The integration strategies developed to solve a scalar complex
second order ODE (\ref{code}), require a two parameter complex group (see,
e.g., \cite{olver,ib}) called $\textbf{G}_{2}$. The integrable forms of
complex second order equations admitting
$\textbf{G}_{2}$ are given in the following Table I.\\
\begin{center}
Table I: Lie canonical forms of complex scalar equations
\begin{tabular}{c c c} \hline\hline
Type & complex symmetry generators & Representative equations \\
\hline\hline
$I$ & $Z_{1}=\partial_{x},~Z_{2}=\partial_{u}$, & $u^{\prime\prime}=w(u^{\prime})$ \\
$II$ & $Z_{1}=\partial_{u},~Z_{2}=x\partial_{u}$, & $u^{\prime\prime}=w(x)$ \\
$III$ & $Z_{1}=\partial_{u},~Z_{2}=x\partial_{x}+u\partial_{u}$, & $xu^{\prime\prime}=w(u^{\prime})$ \\
$IV$ & $Z_{1}=\partial_{u},~Z_{2}=u\partial_{u}$, & $u^{\prime\prime}=u^{\prime}w(x)$ \\
\hline
\end{tabular}
\end{center}

The following theorem summarizes the complex method to solve systems
of ODEs due to their correspondence with the complex scalar solvable
second order ODEs.
\newline\newline
\textbf{Theorem.}\quad A system of two second order ODEs
(\ref{odesys}) is
solvable regardless of the number of symmetries if the corresponding complex scalar equation (\ref{code}) is:\\
(i) \emph{integrable, i.e., it has a two parameter group} $\textbf{G}_2$; or \\
(ii) \emph{linearizable via invertible complex point transformations of
the form} $\mathcal{L}_{2}:(x,u)\rightarrow(\chi,U)$.

Now we give an example of a solvable system which is neither complex linearizable nor linearizable. \newline\newline
\textbf{Example:}\quad Consider a nonlinear coupled system
\begin{eqnarray}
f_1^{\prime\prime}=\frac{(f_1^2-f_2^2)f_1^{\prime}}{(f_1^2-f_2^2)^2+4f_1^2f_2^2}+
\frac{2f_1f_2f_2^{\prime}}{(f_1^2-f_2^2)^2+4f_1^2f_2^2},\nonumber\\
f_2^{\prime\prime}=\frac{(f_1^2-f_2^2)f_2^{\prime}}{(f_1^2-f_2^2)^2+4f_1^2f_2^2}-
\frac{2f_1f_2f_1^{\prime}}{(f_1^2-f_2^2)^2+4f_1^2f_2^2},
\label{sys1}
\end{eqnarray}
which has a two$-$dimensional algebra $[X_1,X_2]=2X_1$, where
\begin{eqnarray}
X_1= \frac{\partial}{\partial x},\quad X_2=
2x\frac{\partial}{\partial x}+ f_1\frac{\partial}{\partial
f_1}+f_2\frac{\partial}{\partial f_2}.
\end{eqnarray}
Using standard Lie analysis it is not straight-forward to carry out
integration of this system. Here we highlight the crucial steps
involved in using the complex transformations (\ref{ctran}) in the
form of invariants and differential invariants of symmetries. We
first observe that $\omega_1$ and $\omega_2$ in (\ref{odesys}),
given by the right hand sides of system (\ref{sys1}) satisfy
CR-equations (\ref{cr0}) therefore system (\ref{sys1}) can be mapped
to a scalar complex equation. Indeed, the equation
\begin{eqnarray}
u^{\prime\prime}=\frac{u^{\prime}}{u^2}, \label{code1}
\end{eqnarray}
corresponds to system (\ref{sys1}) and it has two complex symmetries
$X_1$ and $2x\partial_x+u\partial_u$, therefore (\ref{code1}) has a
solvable Lie algebra $\textbf{G}_2$. The integration of above
equation can be obtained using both approaches, canonical
coordinates or differential invariants. Since scaling is inherited
under $X_1$ therefore we employ canonical coordinates relative to
symmetry $X_1$. The canonical transformation
\begin{eqnarray}
\chi=u, \quad \psi=x,\quad
U(\chi)=\frac{d\psi}{dx}=\frac{1}{u^{\prime}} ,\nonumber
\end{eqnarray}
convert (\ref{code1}) into a first-order equation
\begin{eqnarray}
U^{\prime}=\frac{U^2}{\chi^2},
\end{eqnarray}
which upon realification yields a system of partial differential
equations as $\chi$ is a complex independent variable. That is why
the system (\ref{sys1}) is reduced to a pair of first-order partial
differential equations rather ODEs. This is a similar situation that
arises in complex linearization except the difference that here the
target equation is a reduced solvable ODE not a linear equation. By
integrating the above equation and using invertible transformation
we obtain the solution
\begin{eqnarray}
2c_1f_1+\ln\left( (c_1f_1-1)^2+c_2^2f_2^2\right )-2c_1^2x-2c_1^2c_2=0,\nonumber\\
c_1f_2+\arctan{ \left( \frac{c_1f_1-1}{c_1f_2}\right)}=0,
\end{eqnarray}
of system (\ref{sys1}).

\section*{Acknowledgment} The authors are grateful to Fazal M
Mahomed (South Africa) for useful comments and discussion on this work. We thank the referees for their useful comments which have improved the paper.

\end{document}